\documentclass[preprint, prX]{revtex4}

\usepackage{amsmath}    
\usepackage{graphicx}   
\usepackage{verbatim}   
\usepackage{color}      
\usepackage{subfigure}  
\usepackage{hyperref}   
\usepackage{amssymb}    
\usepackage{epsfig}
\usepackage{graphics,graphicx}
\usepackage{setspace}
\usepackage{url}
\usepackage{algorithm,algorithmic}

\begin{document}

\begin{abstract}
Differential evolution was developed for reliable and versatile function optimization. It has also become interesting for other domains because of its ease to use. In this paper, we posed the question of whether differential evolution can also be used by solving of the combinatorial optimization problems, and in particular, for the graph coloring problem. Therefore, a hybrid self-adaptive differential evolution algorithm for graph coloring was proposed that is comparable with the best heuristics for graph coloring today, i.e. Tabucol of Hertz and de Werra and the hybrid evolutionary algorithm of Galinier and Hao. We have focused on the graph 3-coloring. Therefore, the evolutionary algorithm with method SAW of Eiben et al., which achieved excellent results for this kind of graphs, was also incorporated into this study. The extensive experiments show that the differential evolution could become a competitive tool for the solving of graph coloring problem in the future.

\textit{To cite paper as follows: Fister, I., Brest, J.: Using differential evolution for the graph coloring. In: Proceedings of IEEE SSCI2011 Symposium Series on Computational Intelligence, Piscataway, pp. 150--156 (2011).
}

\end{abstract}

\title{Using Differential Evolution for the Graph Coloring}

\author{Iztok Fister}
\altaffiliation{University of Maribor, Faculty of electrical engineering and computer science
Smetanova 17, 2000 Maribor}
\email{iztok.fister@uni-mb.si}

\author{Janez Brest}
\altaffiliation{University of Maribor, Faculty of electrical engineering and computer science
Smetanova 17, 2000 Maribor}
\email{janez.brest@uni-mb.si}

\maketitle

\section{Introduction}
The graph coloring problem (GCP) can be informally defined as: How to color the graph $G$ with the $k$-colors so that none of the vertices connected with an edge is colored with the same color. The minimum number of colors $k$ needed to color the graph $G$ is called the $chromatic\ number\ \chi$~\cite{book:Kubale2004}. The graph is $k$-colorable if it has a $k$-coloring, i.e. it can be colored with $k$ colors. The problem of finding the correct 3-colorability of the graph $G$ is $\mathit{NP}$-complete~\cite{book:Garey1979}, while the problem of finding the chromatic number of graph $G$ is $\mathit{NP}$-hard~\cite{book:Bondy2008}.

To date, many approaches have been proposed to solve the GCP~\cite{art:Galinier2006,art:Malaguti2010}. The most natural approach to coloring the vertices of graph is in a greedy fashion~\cite{book:Bondy2008}. Here, the vertices are ordered into a permutation and colored sequentially. Although the greedy algorithm, in general, does not find an optimal number of colors, Grimmet and McDiarmid~\cite{inproc:Grimmet1975} showed that this number is no more than about twice the chromatic number. In fact, the quality of greedy coloring depends on the permutation of the vertices. To improve the coloring, Brelaz~\cite{art:Brelaz1979} proposed the algorithm DSatur, which utilizes a dynamic ordering of vertices. This ordering depends on the number of distinctly colored vertices adjacent to the vertex $v$, i.e. the \textit{saturation degree} $\rho(v)$.

While exact algorithms can solve instances up to 100 vertices, heuristic methods are developed for larger instances~\cite{art:Galinier2006}. These can be divided into local search methods~\cite{book:Aarts1997} and hybrid algorithms, which combine the local search methods with evolutionary algorithms~\cite{art:Lu2010}. One of the first local search heuristics for the GCP was Tabucol~\cite{art:Hertz1987}, which uses the tabu search algorithm proposed by Glover~\cite{art:Glover1986}. This heuristic approach was later improved by Fleurent and Ferland~\cite{art:Fleurent1996} and Galinier and Hao~\cite{art:Galinier1999}. The hybrid evolutionary algorithm (HEA) of the later authors is still the best algorithm for GCP today, where the Tabucol that acts as a local search heuristic is guided by the evolutionary meta-heuristic.

The differential evolution (DE) is a population based, very effective global optimization algorithm that was proposed by Price and Storn~\cite{book:Price2005} in 1995. In many ways, DE~\cite{art:Neri2009} is a typical evolutionary algorithm with differential mutation, crossover and selection. Usually, DE was applied to continuous function optimization problems. On the other hand, many attempts to solve the combinatorial optimization problems with DE have been tackled in recent years. Pan et al. in~\cite{art:Pan2008} proposed the discrete differential evolution (DDE) algorithm for the permutation flow scheduling problem. The DDE, hybridized with a variable neighborhood search~\cite{art:Mladenovic1997}, was applied to the generalized assignment problem (GAP) by Tasgetiren et al. in 2009~\cite{inproc:Tasgetiren2009}. Recently, Tasgetiren et al. in 2010~\cite{inproc:Tasgetiren2010} proposed an ensemble of DDE, where each individual is assigned to one of the three distinct DE strategies. A comprehensive survey of the main methods in DE can be found in~\cite{art:Suganthan2011}.

In this paper, we try to answer the question of whether the DE could be used to solve the GCP. We are focused on graph 3-coloring (3-GCP), which is a special kind of $k$-GCP, where $k=3$. Moreover, the 3-GCP can be defined as a constraint optimization problem where the real-valued weights assigned to each vertex establish how difficult the particular vertex is to color. The values of the weights determine the order in which the vertices need to be colored. In other words, these defines the permutation of vertices that determines the sequence in which the constraints should be resolved. In line with this, we translate the definition of 3-GCP to real-valued parameter optimization.

Obviously, the proposed DE is used as a generator of new solutions in the continuous search space, whereas a traditional DSatur algorithm was used for the decoding of vertex permutation to a coloring. This approach is not new, because it was already used by Eiben et al.~\cite{art:Eiben1998}. Despite DSatur, the greedy algorithm was employed as a decoder. Furthermore, the self-adaptation of the parameters F and CR, which controls a differential mutation and crossover, was applied. These parameters are changed using the uncorrelated mutation as proposed by Eiben and Smith~\cite{book:Eiben2003}, which has shown good results by solving of particular combinatorial optimization problems~\cite{incol:Smith2008}. Finally, this algorithm was hybridized with a local search heuristic and, therefore, in this paper, it is denoted as HSA-DE. The local search heuristic is novel and represents a kind of swap heuristic~\cite{book:Eiben2003}. Extensive experiments were conducted on a collection of random 3-colorable graphs generated by the Culbersone graph generator~\cite{website:Culberson}. Note that the selected test suite was the same as by Eiben et al. in~\cite{art:Eiben1998}. Unfortunately, the Dimacs Benchmark suite that consists of the graph coloring problems does not contain any 3-colorable graph.

The structure of the rest of paper is as follows: In Section 2, the 3-GCP is discussed in details. The HSA-DE is described in Section 3, while the experiments and results are presented in Section 4. Through experiments, the phase transition phenomenon is exposed, which refers to the "easy-hard-easy" transition regions, where the problem goes from the status "easy-to-solve" to "hard-to-solve" and conversely~\cite{art:Turner1988}. The paper is concluded with a discussion about the quality of results and directions for the further work.




\section{Graph 3-coloring}

Graph 3-coloring of an undirect graph $G=(V,E)$, where $V$ denotes a finite set of vertices and $E$ a finite set of unordered pairs of vertices $(v_{i},v_{j})$ named edges for $i=1 \ldots n \wedge j=1 \ldots n$, is a mapping $\mathcal{C}:V\rightarrow C$, where $C=\{1,2,3\}$ denotes a set of 3-colors, i.e. the assignment of 3-colors to the vertices of $G$. The proper 3-coloring is obtained when no two adjacent vertices are assigned to the same color~\cite{book:Bondy2008}.

3-GCP can be formally defined as a constraint satisfaction problem (CSP) that is represented as the pair $\langle S,\phi \rangle$, where $S$ denotes the free optimization problem (FOP) and $\phi$ a Boolean function on $S$ (also feasibility condition). For graph 3-coloring the candidate solution $x \in S$ is represented as a permutation of vertices, while the feasibility condition $\phi$ is  expressed as a conjunction of all constraints $\phi (x)= \wedge_{v_{i} \in V} C^{v_{i}}(x)$, where $C^{v_{i}}$ is defined as the set of constraints including the vertex $v_{i}$.

Typically, constraints are handled in DE indirectly in the sense of the penalty function that transforms the CSP into FOP~\cite{book:Eiben2003}. That is, the infeasible solutions that are far away from a feasible region are specified by higher penalties. The penalty function that is used also as objective function here, is expressed as:

\begin{equation}
\label{eq:penalty}
 f(x)=\sum_{i=0}^{n} \psi(x,C^{v_{i}}),
\end{equation}

\noindent where the function $\psi(x,C^{v_{i}})$ is defined as:

\begin{equation}
\label{eq:viol}
\psi(x,C^{v_{i}})=\left\{\begin{matrix}
1 & \textup{if}\ x\ \textup{violates\ at\ least\ one\ } c_{j}\in C^{v_{i}}, \\
0 & \textup{otherwise}.
\end{matrix}\right.
\end{equation}

In fact,~(\ref{eq:penalty}) can be used as a feasibility condition in the sense that $\phi(x)=true$ if and only if $f(x)=0$. Note that this equation evaluates the number of constraint violations and determines the quality of solution $x \in S$ in a permutation search space that is huge, namely $|S|=n!$.

\section{Hybrid Self-adaptive Differential Evolution for Graph 3-coloring}

The DE is simple, powerful, population-based direct search algorithm for the globally optimization of functions~\cite {book:Price2005}. Like other population-based methods, the DE starts with the population of randomly generated initial points in search space. Candidate solutions are modified with the scaled difference of two randomly selected solutions. A trial solution is generated by adding of this difference to a third, randomly selected solution. This trial solution competes with the corresponding solution for survivor. In fact, the solution with the lower value of objective function survives and can evolve in the next generation.

The candidate solution of the DE for 3-GCP is represented as a real-valued vector $Y_{i}^{(t)}=\{w_{i,j}^{(t)}\}$ for $i=1 \ldots \mathit{NP} \wedge j=1 \ldots n \wedge t=1 \ldots \mathit{FEs_{max}}$, where $\mathit{NP}$ denotes the population size, $n$ the number of vector components and $\mathit{FEs_{max}}$ the maximal number of objective function evaluations. The vector components $w_{j}^{(t)} \in [b_{j,L},b_{j,U}]$ representing weights are limited by $n$-dimensional vectors $b_{L}$ and $b_{U}$ that indicate the lower and upper bounds, respectively. Note that the weights determine how hard a corresponding vertex is to color and, thus, define the order in which the vertices will be colored. Each weight is initialized as follows (DE/rand/1/bin):

\begin{equation}
\label{eq:init}
 w_{i,j}^{(0)}=rand(0,1) \cdotp (b_{j,U}-b_{j,L})+b_{j,L},\ \ \ \textnormal{for}\ j=1 \ldots n.
\end{equation}

The DE supports a differential mutation, a differential crossover and a differential selection. In particular, the differential mutation randomly selects two solutions and adds a scaled difference between these to the third solution. This mutation can be expressed as follows:

\begin{equation}
\label{eq:de_mut}
 u_{i}^{(t)}=w_{r0}^{(t)}+F \cdotp (w_{r1}^{(t)}-w_{r2}^{(t)}),\ \ \ \textnormal{for}\ i=1 \ldots NP,
\end{equation}

\noindent where $F \in [0.1,1.0]$ denotes the scaling factor as a positive real number
that scales the rate of modification while $r0,\ r1,\ r2$ are randomly selected vectors in the interval $1 \ldots NP$.

Uniform crossover is employed as a differential crossover by the DE. The trial vector is built out of parameter values that have been copied from two different solutions.  Mathematically, this crossover can be expressed as follows:

\begin{equation}
\label{eq:de_xover}
 z_{i,j}=\begin{cases}
          u_{i,j}^{(t)} & \textnormal{rand}_{j}(0,1) \leq \mathit{CR} \vee j=j_{rand}, \\
		  w_{i,j}^{(t)} & \text{otherwise} ,
        \end{cases}
\end{equation}

\noindent where $\mathit{CR} \in [0.0,1.0]$ controls the fraction of parameters that are copied to the trial solution. Note, the relation $j=j_{rand}$ assures that the trial vector is different from the original solution $Y^{(t)}$.

Mathematically, differential selection can be expressed as follows:

\begin{equation}
\label{eq:de_sel}
 w_{i}^{(t+1)}=\begin{cases}
          z_{i}^{(t)} &\text{if } f(Z^{(t)}) \leq f(Y_{i}^{(t)}), \\
		  w_{i}^{(t)} &\text{otherwise}\,.
        \end{cases}
\end{equation}

Note that the trial vector $Z^{(t)}$ replaces the target vector $Y_{i}^{(t)}$ when the objective function value of the trial vector $f(Z^{(t)})$ is equal or lower than the objective function value of the target vector $f(Y_{i}^{(t)})$. The objective function is calculated according to~(\ref{eq:penalty}).

\subsection{Self-adaptation}

In general, the behavior of the DE depends on the control parameters $F$ and $\mathit{CR}$. However, setting these parameters depends on a problem that is solved by the DE. On the other hand, the initial setting of these parameters can become poor in matured generations. Therefore, the self-adaptation of the control parameters~\cite{book:Baeck1996,art:Brest2006} is introduced into the DE for 3-GCP.

The self-adaptive DE for 3-GCP represents an extension of the classical DE algorithm. Thus, a representation of the solution is expanded by four bytes, i.e. the solution is now represented as a tuple $Y_{i}^{(t)}= \langle w_{i,1}^{(t)}, \ldots , w_{i,n}^{(t)},F_{i}^{(t)}, \sigma_{i,0}^{(t)}, \mathit{CR}_{i}^{(t)}, \sigma_{i,1}^{(t)} \rangle$. The additional elements represent: the control parameter $F_{i}^{(t)}$ with the corresponding mutation strength $\sigma_{i,0}^{(t)}$ and the control parameter $\mathit{CR}_{i}^{(t)}$ with its mutation strength $\sigma_{i,1}^{(t)}$. Furthermore, the parameter $F_{i}^{(t)}$ is mutated according to~(\ref{eq:sigma0})~\cite{book:Eiben2003}:

\begin{equation}
\label{eq:sigma0}
\begin{split}
 \sigma_{i,0}^{(t+1)}=\sigma_{i,0}^{(t)} \cdotp exp(\tau^{'} \cdotp N(0,1)+\tau \cdotp N_{i}(0,1)), &\\
 F_{i}^{(t+1)}=F_{i}^{(t)}+ \sigma_{i,0}^{(t+1)} \cdotp N_{i}(0,1)), &
\end{split}
\end{equation}

\noindent while the parameter $\mathit{CR}_{i}^{(t)}$ is mutated according to~(\ref{eq:sigma1}):

\begin{equation}
\label{eq:sigma1}
\begin{split}
 \sigma_{i,1}^{(t+1)}=\sigma_{i,1}^{(t)} \cdotp exp(\tau^{'} \cdotp N(0,1)+\tau \cdotp N_{i}(0,1)), & \\
 \mathit{CR}_{i}^{(t+1)}=\mathit{CR}_{i}^{(t)}+ \sigma_{i,1}^{(t+1)} \cdotp N_{i}(0,1)). &
\end{split}
\end{equation}

\indent The parameters $\tau^{'} \propto 1/\sqrt{2 \cdotp n}$ and $\tau \propto 1/ \sqrt{2 \cdotp \sqrt{n)}}$ denote the learning rates. To keep the mutation strengths $\sigma_{i,0}^{(t+1)}$ and $\sigma_{i,1}^{(t+1)}$ greater than zero, the following equations are used:

\begin{equation}
\label{eq:eps}
\begin{split}
 \sigma_{i,0}^{(t+1)}<\epsilon_{0} \Rightarrow \sigma_{i,0}^{(t+1)}=\epsilon_{0}\ & and  \\
 \sigma_{i,1}^{(t+1)}<\epsilon_{0} \Rightarrow \sigma_{i,1}^{(t+1)}=\epsilon_{0}. &
\end{split}
\end{equation}

The parameters $\sigma_{i,0}^{(t)}$ and $\sigma_{i,1}^{(t)}$ determines a region in which an evolutionary search can progress and, therefore, a correct setting of their starting values $\sigma_{i,0}^{(0)}$ and $\sigma_{i,1}^{(0)}$ have the significant influence on the performance of this algorithm.

\subsection{Genotype-Phenotype Mapping}

Conventional DE algorithms explore continuous search space (also known as genotype space). However, 3-GCP is a combinatorial problem and has a discrete problem space (also known as phenotype space). In contrast to many graph coloring algorithms that encode solutions as binary strings~\cite{art:Greenwood2009}, the proposed self-adaptive DE algorithm employs a special genotype-phenotype mapping. This acts as follows: A solution of the 3-GCP is defined within the self-adaptive DE algorithm as a tuple $Y_{i}^{(t)}= \langle w_{i,1}^{(t)}, \ldots , w_{i,n}^{(t)},F_{i}^{(t)}, \sigma_{i,0}^{(t)}, \mathit{CR}_{i}^{(t)}, \sigma_{i,1}^{(t)} \rangle$ for $i=1 \ldots \mathit{NP}$. That is, the solution is described by this DE algorithm indirectly. All candidate solutions thereby define the genotype search space. On the other hand, a quality of these solutions is evaluated in the phenotype problem space.

The genotype-phenotype mapping consists of two phases:
\begin{itemize}
  \item transformation of weights in permutation of vertices,
  \item decoding of solution by DSatur construction heuristic,
\end{itemize}
In the first phase, only weights are taken from the solution $Y_{i}^{(t)}$. Then, these are put in descending order and, according to their ordering, a permutation of vertices $X_{i}^{(t)}=\{v_{i,j}^{(t)}\}$ for $j=1 \ldots n$ is constructed. From this permutation, the DSatur construction heuristic decodes a coloring $C_{i}^{(t)}=\{c_{i,j}^{(t)}\}$ for $i=1 \ldots \mathit{NP} \wedge j=1 \ldots n$, where $c_{i,j}^{(t)} \in \{ 1,2,3\}$ represents the solution in a phenotype space. Note that the phenotype space is much smaller than the genotype space, namely $|C|=3^{n}$.

After decoding, the quality of the solution is evaluated by the counting the number of constraint violations according to~(\ref{eq:penalty}).

\subsection{Hybridization with Local Search}

The current solution of the self-adaptive DE can be improved with a local search heuristic. The local search heuristic tries to improve the current solution with small adjustments to weights until improvements are perceived. The proposed algorithm was hybridized by the ordering by saturation heuristic. In line with hybridization, the proposed algorithm was named: the hybrid self-adaptive differential evolution algorithm (HSA-DE).

\begin{figure*}[!htb]
    \begin{center}
        \includegraphics[scale=0.72]{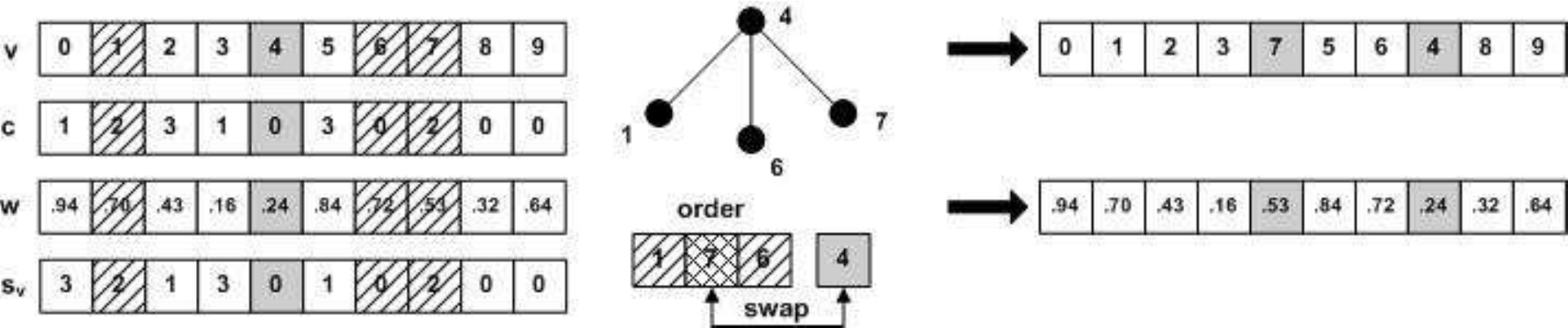}
        \caption{Order local search heuristic.}
        \label{pic:slika_1}
    \end{center}
\end{figure*}

The proposed heuristic acts as follows: Initially, it locates the first uncolored vertex. Then, a set of neighbors to this vertex is determined. The set that represents a neighborhood are put in descending ordered with regard to the values of saturation degrees. In the end, the neighbor with the highest value of saturation degree is swapped with the uncolored vertex. Actually, the neighborhood is defined by a set of possible 1-$swap$ strokes.

The way, in which this heuristic operates is illustrated in Fig.~\ref{pic:slika_1}, where $v$ denotes the permutation of weights, $c$ the corresponding coloring, $w$ the weights and $d$ the saturation degrees. Here, the first uncolored vertex 4 is shadowed, while its adjacent vertices 1, 6 and 7 are hatched. Because the vertices 1 and 7 have the same saturation degree, the vertex 7 is selected randomly. Finally, the vertices 4 in 7 are swapped (right-hand side of Fig.~\ref{pic:slika_1}).

\section{Experiments and results}

The goal of an experimental work was to show that the HSA-DE can be used to solve the 3-GCP. Therefore, a comparison with algorithms:
\begin{itemize}
  \item Tabucol~\cite{art:Hertz1987},
  \item the hybrid evolutionary algorithm (HEA)~\cite{art:Galinier1999},
  \item the evolutionary algorithm with the method SAW (EA-SAW)~\cite{art:Eiben1998} and
  \item the hybrid self-adaptive differential evolution algorithm (HSA-DE),
\end{itemize}
\noindent was accomplished. The first two algorithms are some of the best known algorithms for the $k$-GCP today. Furthermore, the implementation of both are published on internet sites~\cite{website:Chiarandini} as help to developers of new graph coloring algorithms. Interestingly, the implementation of EA-SAW can be found online~\cite{website:Hemert} as well. However, the later algorithm was designed especially for the 3-GCP. In our opinion, the paper of Eiben et al.~\cite{art:Eiben1998}, in which it was described, is one of the most complete publications about 3-GCP today. Thus, this paper was taken as a reference for comparison.

The characteristics of the proposed HSA-DE, as used in the experiments, are as follows. A population size was set to 15 because this value represents a good selection as was indicated during the experimental work. In general, the strategy DE/rand/1/bin was employed. The local search heuristic was applied with the probability $p_{ls}=0.02$. As previously mentioned, the most attention was given to setting the correct starting values of mutation strengths $\sigma_{i,0}^{(0)}$ and $\sigma_{i,1}^{(0)}$. After preliminary experimentation, the values $\sigma_{0}^{(0)}=30.0$ and $\sigma_{1}^{(0)}=30.0$ were selected. The number of objective function evaluations $\mathit{FEs_{max}}=300,000$ was used as the termination condition. Because of the stochastic nature of the mentioned algorithms, multiple runs were observed, and the number of runs was limited to 25.

Algorithms were compared according to the measures:
\begin{itemize}
  \item success rate (SR) and
  \item average number of objective function evaluations to solution (AES),
\end{itemize}
\noindent where the first measure expresses the ratio of the number of successful runs among all runs and the second reflects the efficiency of the particular algorithm.

\subsection{Test Problems}

A test suite consisted of random graphs generated by Culbersone graph generator~\cite{website:Culberson}. The generated graphs are denoted as $G_{t,n,p,s}$, where $t$ represents the \textit{type of graphs}, $n$ the \textit{number of vertices}, $p$ the probability that two verices $v_{i}$ and $v_{j}$ are connected with an edge $(v_{i},v_{i})$ and $s$ the \textit{seed} of the random graph generator. Indirectly, the parameter $p$ controls the \textit{edge density}. Among the many types of graphs equi-partite, uniform and flat graphs were generated because of the comparison with the selected reference~\cite{art:Eiben1998}. In equi-partite graphs, the vertices are split into 3-color sets and then edges are selected so that the sizes of all 3-color sets are as nearly as possible. In uniform graphs, the vertices are assigned to one of the 3-color sets uniformly and independently, according to the predefined probability~\cite{incol:Chiarandini2010}. The flat graphs are similar to the equi-partite, with an additional property that variations of expected vertex degrees are minimized. All graphs were generated with the number of vertices set at $n=1,000$. These kind of graphs were also referred to as large-scale in~\cite{art:Eiben1998}. The seed $s$ of the random graph generator does not have a notable impact on the results, as noted by Eiben et al.~\cite{art:Eiben1998}. To keep the comparison fairly, the $s=5$ was used.

The most combinatorial optimization problems are sensitive on the phase transition phenomenon that determines the regions, where the problem passes from the state of "solvable" to the state of "unsolvable" and vice versa~\cite{art:Turner1988}. Typically, these regions are connected with some parameter addressing the problem. For the 3-GCP, this parameter is the edge probability $p$. Many authors identified different critical values of this parameter, which determines the phase transition region. For example, Petford and Welsh~\cite{art:Petford1989} stated that this phenomenon occurs when $2pn/3 \approx 16/3$, Cheeseman et al.~\cite{inproc:Cheeseman1991} when $2m/n \approx 5.4$, Eiben et al.~\cite{art:Eiben1998} when $7/n \leq p \leq 8/n$ and Hayes~\cite{art:Hayes2003} when $m/n \approx 2.35$. Note that the parameter $m$ in the above-mentioned formulas denotes the number of edges. In this paper, however, we consider that the hard problems occur when $p=0.007$.

\subsection{Influence of the Edge Density}

In this test, the phenomenon of phase transition was investigated. Therefore, random graphs were generated in which edge probabilities varied from $p=0.004$ to $p=0.014$ in steps of $0.001$. As a result, 11 instances of randomly generated graphs were obtained for each type of the observed graphs.

\begin{figure*}[hbt]	
\centering
\subfigure[SR by equi-partite graphs] {\includegraphics[width=7.2cm]{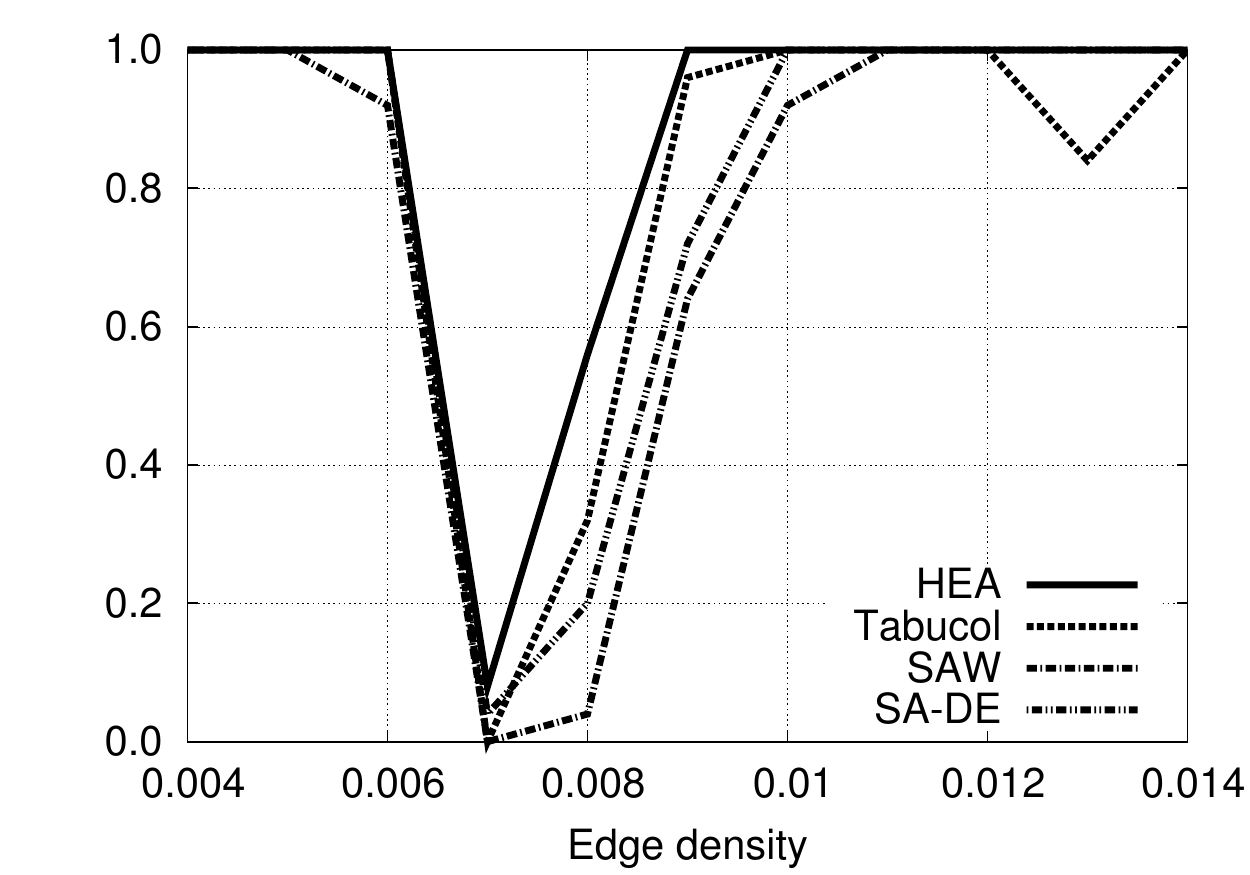}}
\subfigure[AES by equi-partite graphs] {\includegraphics[width=7.2cm]{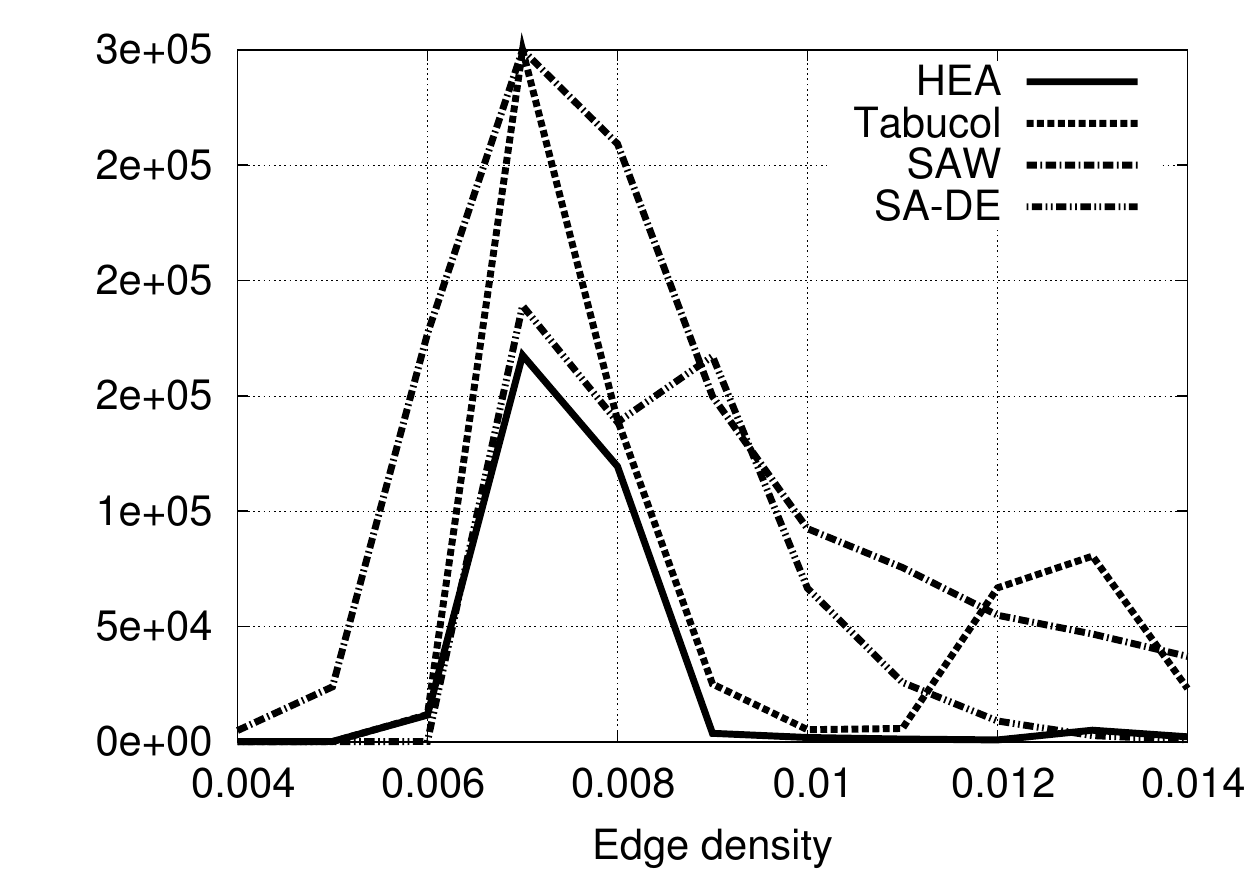}}
\subfigure[SR by uniform graphs] {\includegraphics[width=7.2cm]{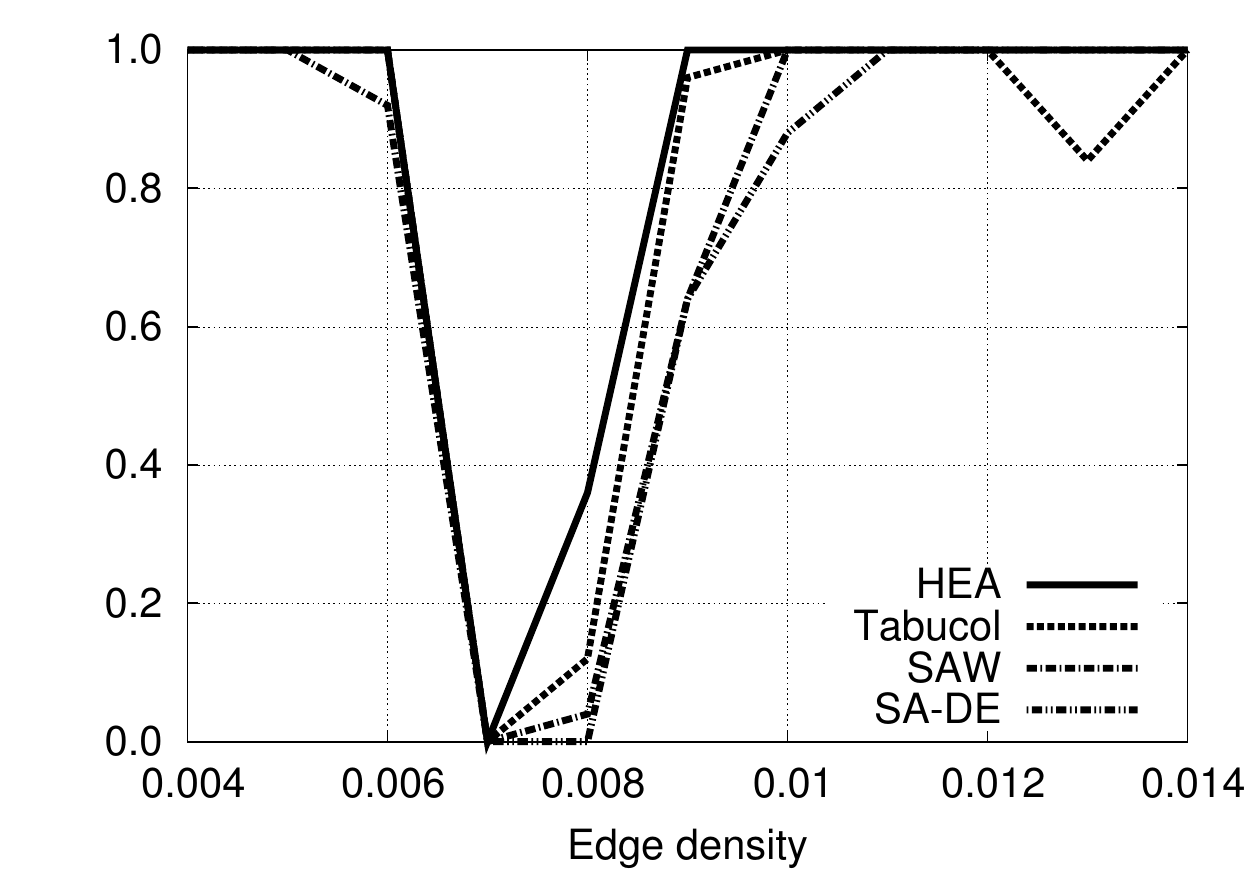}}
\subfigure[AES by uniform graphs] {\includegraphics[width=7.2cm]{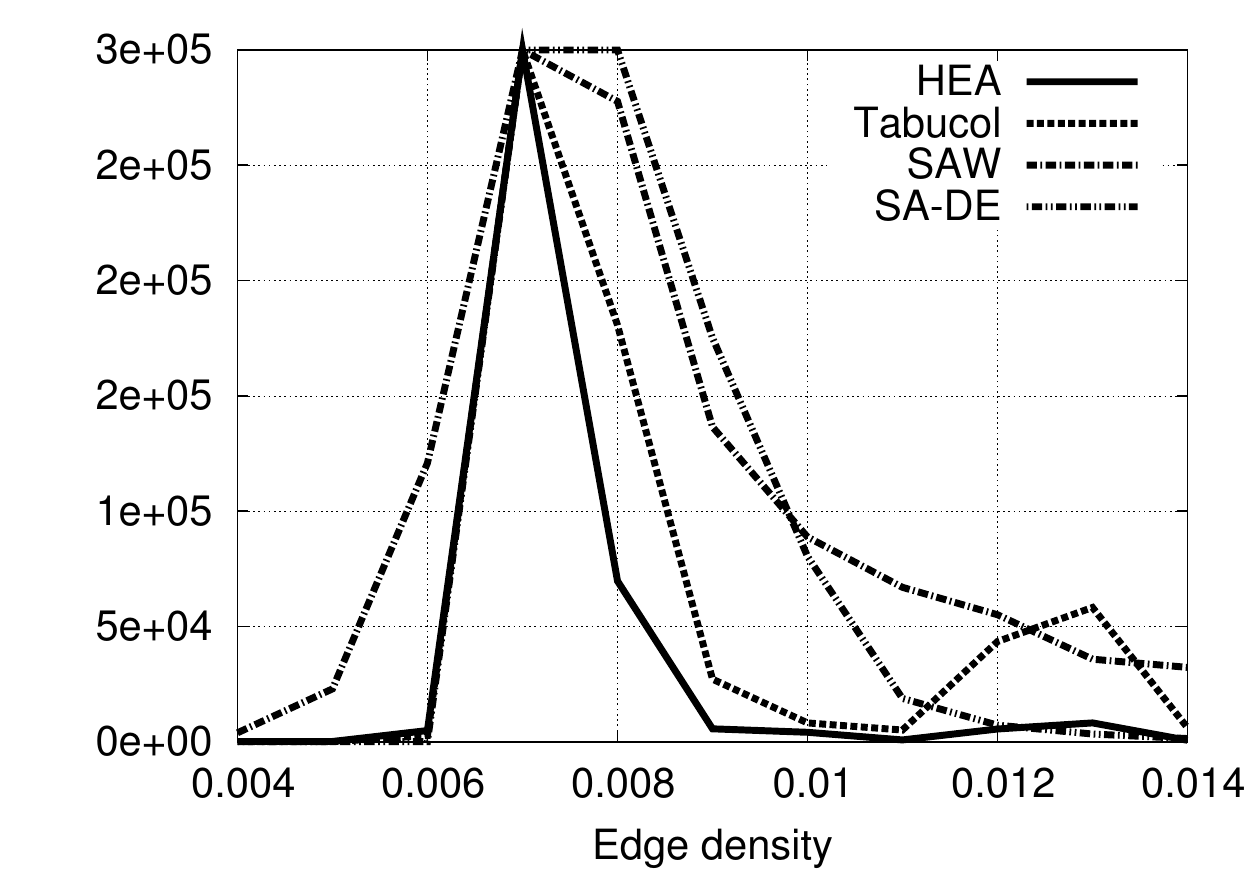}}
\subfigure[SR by flat graphs] {\includegraphics[width=7.2cm]{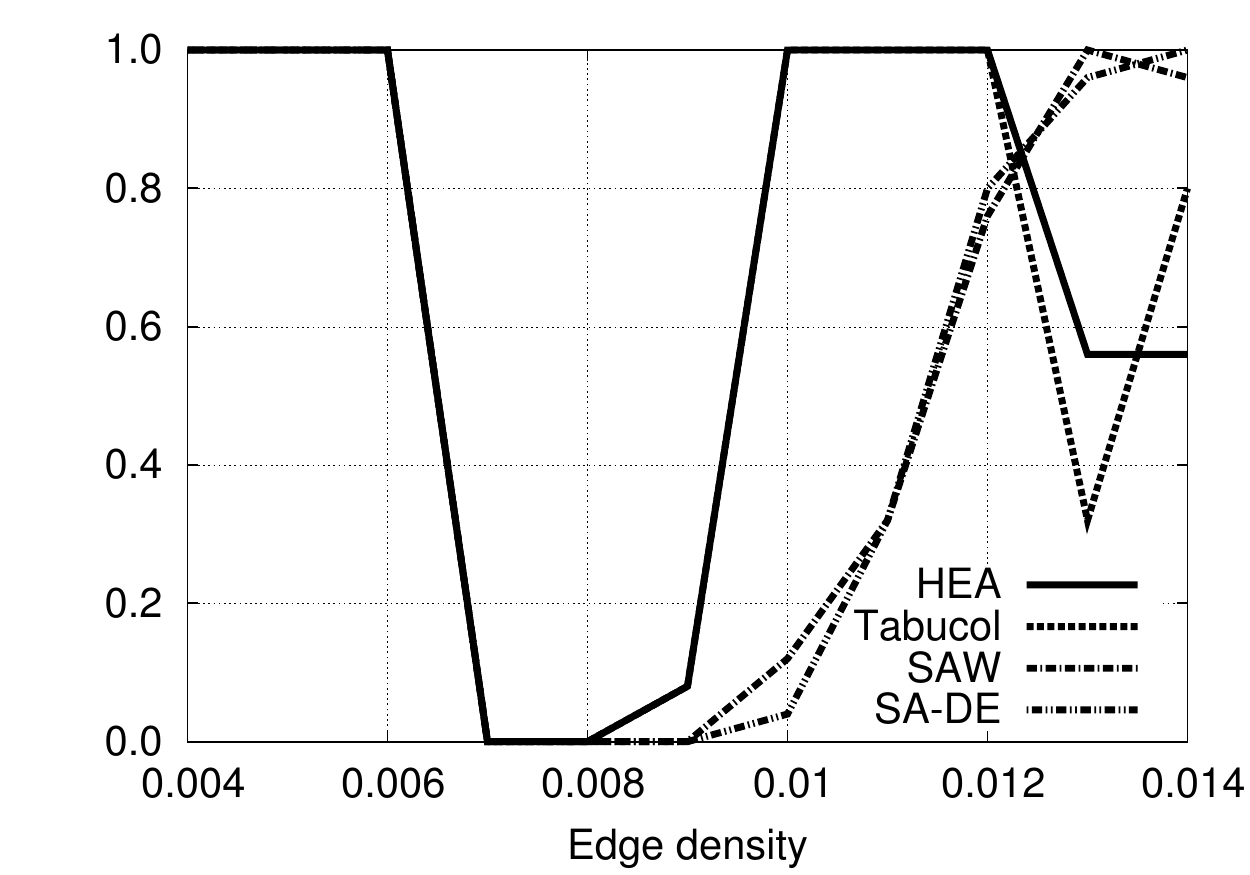}}
\subfigure[AES by flat graphs] {\includegraphics[width=7.2cm]{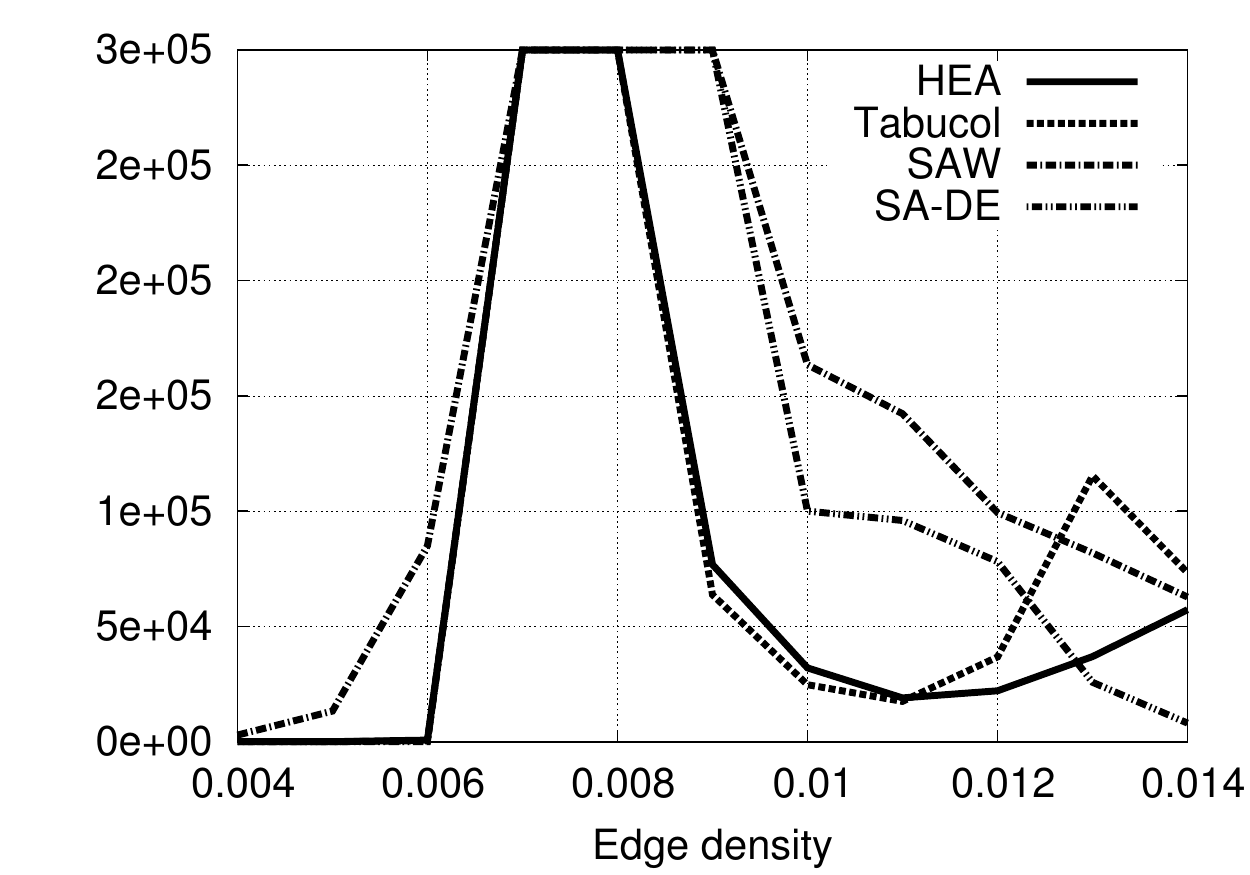}}
\caption{Results of algorithms for 3-GCP solving equi-partite, uniform and flat random graphs}
\label{fig:Sub_1}
\end{figure*}

The results of the comparison are illustrated in Fig.~\ref{fig:Sub_1}. This figure is divided into six graphs, according to the type and measures SR and AES. The graphs capture the results of all 11 instances that were obtained by varying the edge probability $p$.

As can be seen in Fig.~\ref{fig:Sub_1}.a, the HEA is certainly the best algorithm for GCP, based on the measure SR. Although it has some difficulties solving graphs with the edge probability $p=0.007$, where the SR=0.08 was reached and the $p=0.008$, where the SR=0.56 was reached, other graph instances were solved excellently (SR=1). Interestingly, the HSA-DE also solved all graph instances but results by 3-coloring outside of the phase transition region were worse than with Tabucol and better than with EA-SAW. Tabucol and EA-SAW did not solve the graph with $p=0.007$. Moreover, Tabucol also did not solve the graph with $p=0.013$ in all runs that needed to be outside of the phase transition region. This problem is connected with the fact that here the second phase transition exists~\cite{art:Mulet2002} as well. The results, according to the measure AES in Fig.~\ref{fig:Sub_1}.b, show that the HSA-DE, more or less, follows the results of the HEA. Interestingly, results far away from the phase transition region are even better than with the other observed algorithms.

As illustrated by Fig.~\ref{fig:Sub_1}.c, the results of the HSA-DE on uniform graphs are comparable with the results of the EA-SAW, but worse than results of the HEA and Tabucol. Similar to the equi-partite graphs, Tabucol was sensitive in the second phase transition region. According to the AES (Fig.~\ref{fig:Sub_1}.d) the results of the HSA-DE were worse than the results of all other observed algorithms in the vicinity of the phase transition region, but overcame the results of the other algorithms in the region outside of it.

Finally, the HSA-DE come nearer to the results of the EA-SAW, according to the SR, by 3-coloring of flat graphs (Fig.~\ref{fig:Sub_1}.e). However, the results of both algorithms were worse than the results of the algorithms for $k$-GCP. As can be seen in the Fig.~\ref{fig:Sub_1}.e, these types of graphs are the hardest for all of the observed algorithms. Moreover, both algorithms for $k$-GCP are sensitive to the second phase transition as well. According to the AES (Fig.~\ref{fig:Sub_1}), the HSA-DE improved the results of the EA-SAW before and after the phase transition region.

\begin{figure*}[hbt]	
\centering
\subfigure[Number of uncolored vertices by successfully run] {\includegraphics[width=7.2cm]{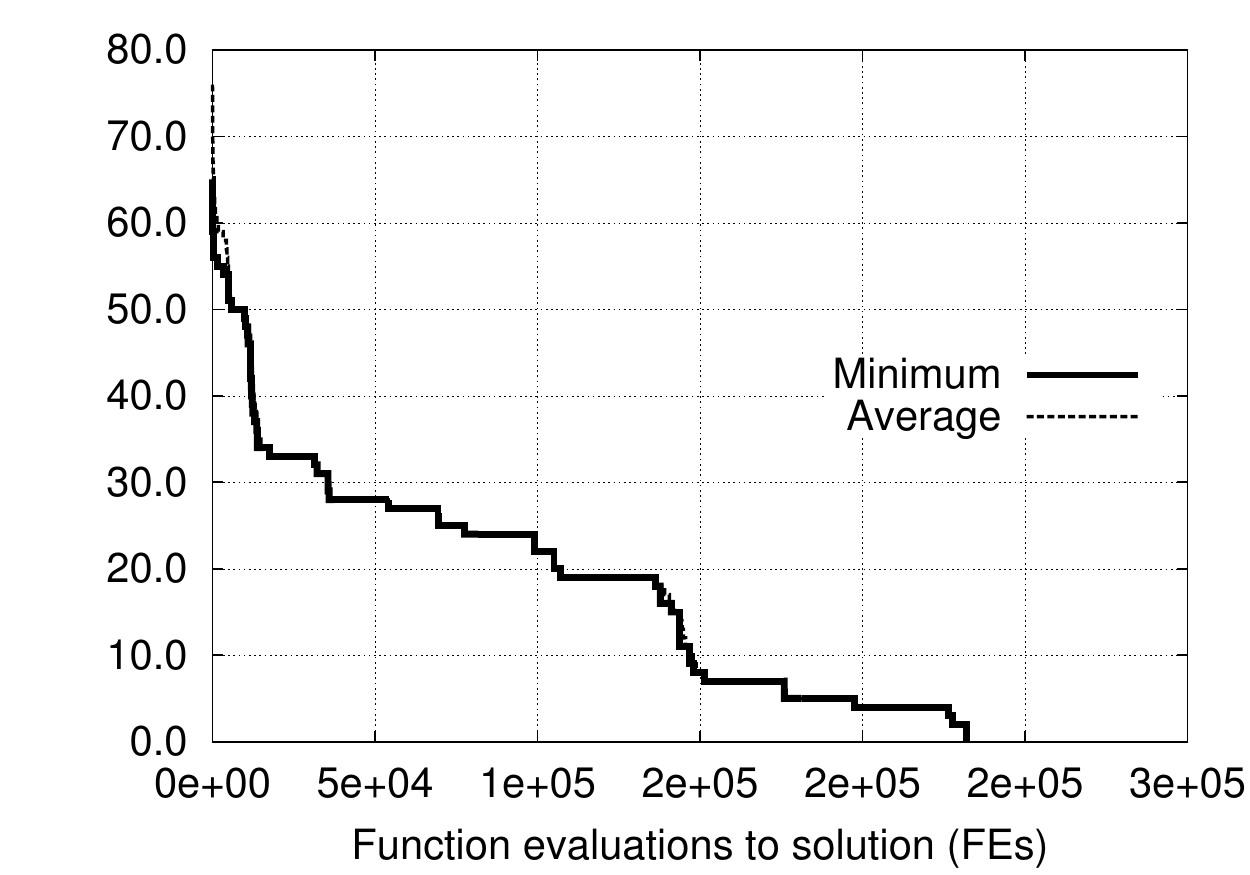}}
\subfigure[Number of uncolored vertices by unsuccessfully run] {\includegraphics[width=7.2cm]{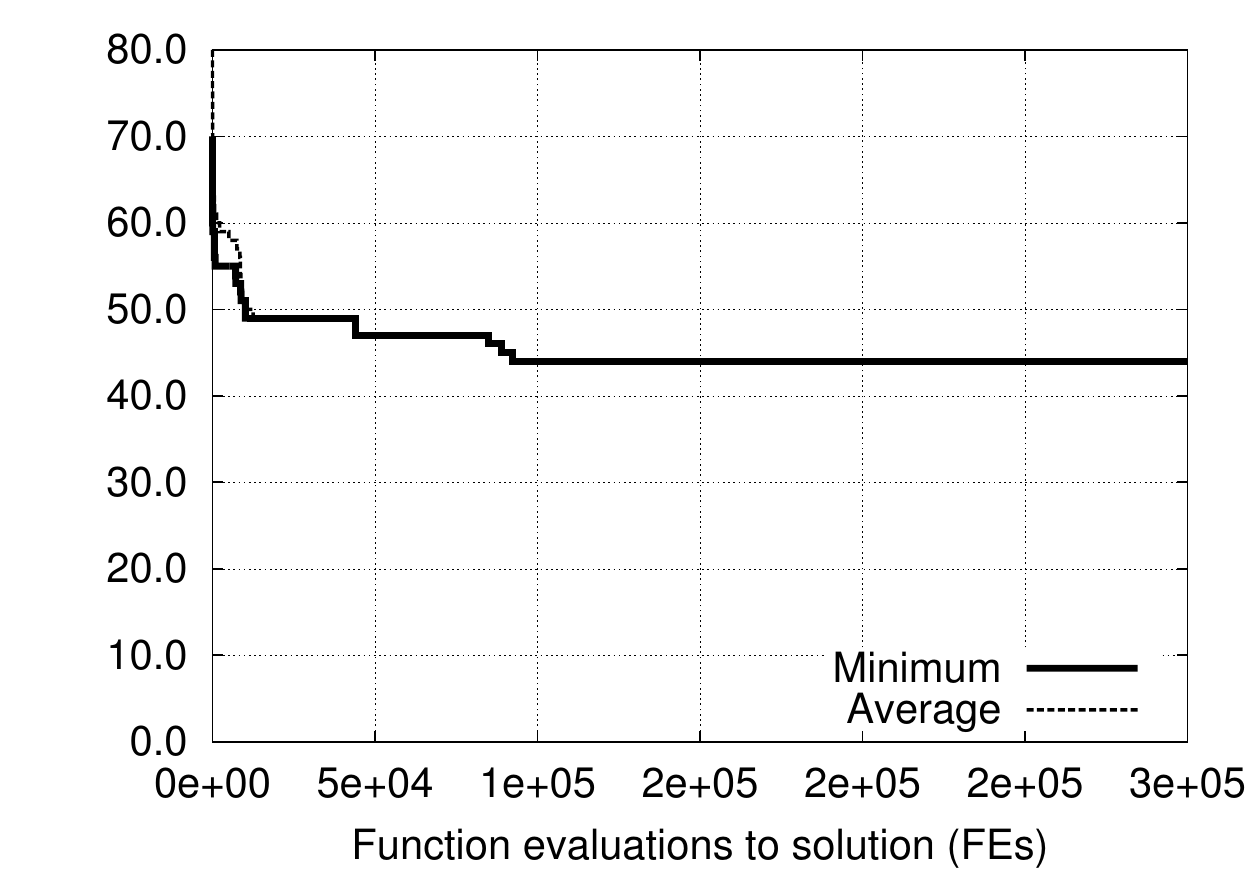}}
\caption{Two independent runs of the HSA-DE}
\label{fig:Sub_2}
\end{figure*}

\subsection{Influence of the Local Search Heuristic}

A goal of this experiment was to investigate the influence of the local search heuristic on results of the HSA-DE. In this sense, the self-adaptive DE without a local search was compared with the original HSA-DE. Because we focused on the behavior of the HSA-DE in the phase transition region, the edge density of graphs were varied from $p=0.006$ to $p=0.009$ with steps of 0.001. Note that the equi-partite graphs were investigated only. The results of the experiment are presented in Table~\ref{tab:tab_1}.

\begin{table}[htb]        
\caption{Influence of the local search by HSA-DE on equi-partite graphs.}
\label{tab:tab_1}
\begin{center}
\begin{tabular}{  l  l  l  l  l  }
\hline
 & HSA-DE & NONE & HSA-DE & LS  \\
p & SR & AES & SR & AES  \\
\hline
0.006 & 1 & 52.76 & 1 & 70.92   \\
0.007 & 0 & 0 & 0.04 & 189,087  \\
0.008 & 0.12 & 267,364 & 0.2 & 138,542  \\
0.009 & 0.44 & 170,447 & 0.72 & 167,178  \\
\hline
avg & 0.39 & 184,465.94 & 0.49 & 123,719.48 \\
\hline
\end{tabular}
\end{center}
\end{table}

From Table~\ref{tab:tab_1} it can be seen that hybridization with a local search was successful when the probability $p_{ls}=0.02$ was used, but increasing of this parameter can lead to worse results, as indicated by experiments.

\subsection{Discussion}

To discover how the HSA-DE explores a search space a deeper analysis of how it acts needs to be performed. Therefore, two independent runs of the HSA-DE were examined by solving the equi-partite graphs with $p=0.008$, where the solution is either found as illustrated in Fig.~\ref{fig:Sub_2}.a or not found as in Fig.~\ref{fig:Sub_2}.b. In these figures, the average number of uncolored vertices in the population and the current minimum number of uncolored vertices are recorded throughout the runs.

In Fig.~\ref{fig:Sub_2}.a, it can be seen that the curve, which denotes the minimum number of uncolored vertices, converges to the optimal value very fast at beginning. Then, the decline is slower. When the solution is not found (Fig.~\ref{fig:Sub_2}.b) the current minimum number of uncolored vertices falls from 80 to 50 very rapidly, but when the evolutionary search progresses, no further improvements are detected. That is: the HSA-DE sinks into the local optima. If we look at the average number of uncolored vertices in the population, we can observe that this value is very close to the best number of uncolored vertices. Because this value can be taken as a measure of population diversity \cite{book:Eiben2003} we can conclude that the population looses its diversity very quickly. As noted by B{\"a}ck in~\cite{book:Baeck1996}, the diversity of the population is a prerequisite for the self-adaptation. Therefore, it needs to be increased. Obviously, diversity can be increased in many ways. For example, it can be done by increasing of the population size $\mathit{NP}$, or conducting a new vectors in the population during the evolutionary process, etc. Although the diversity of the population was increased by applying either of the proposed mechanisms, at this moment, neither of these bring any visible improvement to the results. Therefore, to find the mechanism that would be able to balance an exploration (population diversity) and an exploitation (selection pressure)~\cite{book:Eiben2003} in DE remains the great challenge for future work.

\section{Conclusion}

In summary, the gained results show that the HSA-DE for 3-GCP is better than the EA-SAW by solving of equi-partite graphs and it is comparable by solving of other three graph types. However, both of the mentioned algorithms achieved worse results than algorithms for $k$-GCP, i.e. HEA and Tabucol, by solving of all other three types of observed graphs. Obviously, our intention in this paper was not to improve the results of these algorithms but to show that DE can be applied to the GCP. In this first step, the HSA-DE obtained results that are comparable with EA-SAW. However, in order to develop an algorithm that would be comparable with the best algorithms for $k$-GCP, more investigations would need to be conducted.

\bigskip{\small \smallskip\noindent Updated 30 November 2012.}
\end{document}